\documentclass[a4paper, 12pt]{article}
\usepackage{amssymb}
\usepackage{amsmath}
\usepackage{amscd}
\usepackage{latexsym}
\usepackage{theorem}

\usepackage{color}

\theoremheaderfont{\scshape}
\newtheorem{thm}{\bfseries Theorem}[section]
\newtheorem{prop}[thm]{\bfseries Proposition}
\newtheorem{lemma}[thm]{\bfseries Lemma}
\newtheorem{cor}[thm]{\bfseries Corollary}

\theorembodyfont{\rmfamily}
\newtheorem{defn}[thm]{\bfseries Definition}

\newtheorem{rem}[thm]{Remark}
\newtheorem{notation}[thm]{Notation}
\newtheorem{pf}{Proof.}

\def\m{\mathfrak m}

\def\Hom{\mathrm{Hom}}
\def\pHom{\underline{\mathrm{Hom}}}
\def\Ext{\mathrm{Ext}}

\def\End{\mathrm{End}}
\def\pEnd{\underline{\mathrm{End}}}

\def\ker{\mathrm{Ker}}

\def\coker{\mathrm{Coker}}

\def\depth{\mathrm{depth\,}}

\def\ind{\mathrm{Ind}}
\def\rad{\mathrm{rad}}
\def\G{\mathcal G}
\def\H{\mathcal H}
\def\C{\mathcal C}
\def\A{\mathcal A}
\def\pG{\underline{\mathcal G}}
\def\pH{\underline{\mathcal H}}
\def\pC{\underline{\mathcal C}}
\def\modR{\mathrm{mod}{\hspace{1pt}}R}
\def\modpG{\mathrm{mod}{\hspace{1pt}}\underline{\mathcal G}}
\def\modpH{\mathrm{mod}{\hspace{1pt}}\underline{\mathcal H}}
\def\pmodR{\underline{\mathrm{mod}{\hspace{1pt}}R}}

\def\modC{\mathrm{mod}{\hspace{1pt}}{\mathcal C}}

\def\modpC{\mathrm{mod}{\hspace{1pt}}\underline{\mathcal C}}
\def\ModpC{\mathrm{Mod}{\hspace{1pt}}\underline{\mathcal C}}
\def\modA{\mathrm{mod}{\hspace{1pt}}{\mathcal A}}
\def\ModA{\mathrm{Mod}{\hspace{1pt}}{\mathcal A}}
\def\tr{\mathrm{Tr}}

\def\rad{{\mathrm{rad}}}

\def\qed{$\Box$}
\title{\bf A functorial approach to modules of G-dimension zero \\}
\author{Yuji Yoshino \\
Math. Department, Faculty of Science \\
Okayama University, Okayama 700-8530, Japan \\
\vspace{4pt}
\texttt{yoshino@math.okayama-u.ac.jp}\\
}

\date{\empty}

\begin{document}
\maketitle

\begin{abstract}
Let  $R$  be a commutative Noetherian ring and let  $\G$  be the category of modules of G-dimension zero over  $R$.
We denote the associated stable category by  $\pG$.
We show that the functor category $\modpG$ is a Frobenius category and we argue  how this property could characterize $\G$  as a subcategory of $\modR$.

\vspace{6pt}

\noindent
{\bf Math Subject Classification 2000}: 13C  , 13D 
\end{abstract}


\section{Introduction}
In this paper  $R$  always denotes a commutative Noetherian ring, and  $\modR$  is the category of finitely generated $R$-modules. 

In this paper we are interested in the subcategories $\G$ and $\H$  of $\modR$  that are defined as follows: 
First, $\G$ is defined to be the full subcategory of  $\modR$ consisting of all modules  $X \in \modR$  that satisfy 
$$
\Ext _R ^i(X, R) = 0  \quad  \text{and} \quad 
\Ext _R ^i(\tr X, R) = 0 \quad  \text{for any } \  i >0.
$$
We also define $\H$ to be the full subcategory consisting of all modules with the first half of the above conditions, therefore  a module  $X \in \modR$ is an object in  $\H$ if and only if 
$$
\Ext _R ^i(X, R) = 0  \quad  \text{for any } \  i >0.
$$
Note that $\G \subseteq \H$ and that  $\G$  is called the subcategory of modules of G-dimension zero.
See \cite{ABr} for the G-dimension of modules. 
We optimistically expect that the equality  $\G = \H$  holds in many cases.

The main purpose of this paper is to characterize functorially these two subcategories and to get the conditions under which a subcategory  $\C$ of  $\H$  is contained in  $\G$. 

For any subcategory  $\C$ of  $\modR$, we denote by $\pC$  the associated stable category  and denote by  $\modpC$  the category of finitely presented contravariant additive functors from $\pC$ to the category of Abelian groups. 
See \S 2 for the precise definitions for these associated categories. 
As a first result of this paper we shall prove in \S 3  that the functor category $\modpH$  is a quasi-Frobenius category, while $\modpG$  is a Forbenius category. 
See Theorems \ref{quasi-Frobenius} and \ref{Frob}. 

For any subcategory  $\C$  of $\modR$ which satisfy several admissible mild conditions, we expect that Frobenius and quasi-Frobenius property of the category $\modpC$ will be equivalent to that  $\C$ is contained in  $\G$  and $\H$  respectively. 

Let  $R$  be a henselian local ring, and hence the category  $\modR$ is a Krull-Schmidt category. 
Then we shall prove in  Theorem \ref{characterize H} that  a resolving subcategory $\C$  of  $\modR$ is contained in $\H$ if and only if $\modpC$ is a quasi-Frobenius category. 
 In this sense we can characterize the subcategory  $\H$  as the maximum subcategory  $\C$  for which  $\modpC$  is a quasi-Frobenius category.

For the subcategory  $\G$,  a similar functorial characterization is also possible but with an assumption that the Auslander-Reiten conjecture is true.
See Theorem \ref{moduloAR} for the detail.

In the final section \S 5 we shall prove the main Theorem \ref{characterize G} of this paper,  in which we assert that any resolving subcategory of  $\H$ that is of finite type is contained in  $\G$. 
In particular, if  $\H$ itself is of finite type,  then we can deduce the equality  $\G = \H$.

\section{Preliminary and Notation}

Let   $R$  be a commutative Noetherian ring, and  let  $\modR$  be the category of finitely generated $R$-modules as in the introduction. 

When we say $\C$ is a subcategory of $\modR$, we always mean the following:
\begin{itemize}
\item 
$\C$ is essential in $\modR$, 
i.e. if $X \cong Y$  in $\modR$ and if $X \in \C$, then  $Y \in \C$.
\item 
$\C$ is full in $\modR$, 
i.e. $\Hom _{\C} (X, Y) = \Hom _R(X, Y)$  for $X, Y \in \C$.
\item  
$\C$ is additive and additively closed in $\modR$, 
i.e. for any $X, Y \in \modR$,  $X \oplus Y \in \C$  if and only if  $X \in \C$  and  $Y \in \C$.
\item
$\C$ contains all projective modules in $\modR$.
\end{itemize}

The aim of this section is to settle the notation that will be used throughout   this paper and to recall several notion for the categories associated to a given subcategory.

Let  $\C$  be any subcategory of $\modR$. 
At first, we define the associated stable category  $\pC$ as follows: 
\begin{itemize}
\item
The objects of $\pC$ are the same as those of  $\C$.
\item
For  $X, Y \in \pC$, the morphism set is an $R$-module
$$
\pHom _{R}(X, Y) = \Hom _R (X, Y)/P(X, Y),  
$$
 where  $P(X, Y)$  is the $R$-submodule of  $\Hom _R(X, Y)$  consisting of all $R$-homomorphisms which factor through projective modules.  
\end{itemize}

Of course, there is a natural functor  $\C \to \pC$. 
And for an object  $X$ and a morphism  $f$  in $\C$  we denote their images in  $\pC$  under this natural functor by $\underline{X}$ and  $\underline{f}$. 

\begin{defn}
Let  $\C$  be a subcategory of  $\modR$.
For a module $X$  in $\C$, we take a finite presentation by finite projective  modules 
$$
P_1 \to P_0 \to X \to 0, 
$$ 
and define the transpose  $\tr X$ of  $X$   as the cokernel of  $\Hom _R(P_0, R ) \to \Hom _R (P_1, R)$. Similarly for a morphism $f : X \to Y$  in  $\C$, since it induces a morphism between finite presentations 
$$
\begin{CD}
P_1 @>>> P_0 @>>> X @>>> 0 \\
@V{f_1}VV  @V{f_0}VV  @V{f}VV \\
Q_1 @>>> Q_1 @>>> Y @>>> 0,  \\
\end{CD}
$$
 we define the morphism $\tr f : \tr Y \to \tr X$  as the morphism induced by  $\Hom _R (f_1, R)$.
It is easy to see that $\tr X$ and $\tr f$  are uniquely determined as an object and  a morphism in the stable category  $\pC$, and 
it defines well the functor 
$$
\tr :  (\pC) ^{op} \to \pmodR.
$$

For a module  $X \in \C$  its syzygy module $\Omega X$  is defined by the exact sequence
$$
0 \to \Omega X \to P_0 \to X \to 0, 
$$
where  $P_0$  is a projective module. 
It is also easy to see that $\Omega$  defines a functor 
$$
\Omega : \pC \to  \pmodR.
$$
\end{defn}

We are interested in this paper two particular subcategories and their associated stable categories.

\begin{notation}
We denote by $\G$ the subcategory of  $\modR$ consisting of all modules of G-dimension zero, that is,  a module $X \in \modR$  is an object in $\G$  if and only if 
$$
\Ext _R ^i(X, R) = 0  \quad  \text{and} \quad 
\Ext _R ^i(\tr X, R) = 0 \quad  \text{for any } \  i >0.
$$

We also denote by $\H$ the subcategory consisting of all modules with the first half of the above conditions, that is,  a module  $X \in \modR$ is an object in  $\H$ if and only if 
$$
\Ext _R ^i(X, R) = 0  \quad  \text{for any } \  i >0.
$$
\end{notation}

Of course we have  $\G \subseteq \H$.
Since we have no example of modules in  $\H$  that is not in $\G$, 
we conjecture very optimistically that  $\G = \H$.
This is a main point of this paper, and we argue how they are different or similar by characterizing these subcategories by a functorial method.

Note just from the definition that  $\tr$  gives dualities on  $\pG$  and  $\pmodR$, that is, the first and the third vertical arrows in the following diagram are isomorphisms of categories:
$$
\begin{CD}
(\pG)^{op} \quad  @. {\subseteq} \quad @. (\pH)^{op}  \quad @. {\subseteq} \quad @. (\pmodR)^{op} \\
@V{\tr}VV \quad @. @V{\tr}VV  \quad @. @V{\tr}VV  \\ 
\pG  \quad @. {\subseteq}  \quad @. \tr\pH  \quad @. {\subseteq}  \quad @. \pmodR \\
\end{CD}
$$
Here we note that  $\tr \pH$  is the subcategory consisting of all modules $X$ 
satisfying  $\Ext ^i _R (\tr X, R)=0$  for all $i>0$, hence we have the equality$$
\pG = \pH \cap \tr \pH.
$$
Therefore  $\pG = \pH$  is equivalent to that $\tr \pH = \pH$, that is,  $\pH$  is closed under $\tr$.
Note also that  $\pG$  and $\pH$  are closed under the syzygy functor, i.e. $\Omega\pG = \pG$ and  $\Omega \pH \subseteq  \pH$.

For an additive category  $\A$, a contravariant additive functor from  $\A$  to the category $(Ab)$ of abelian groups is referred to as an $\A$-module, 
and a natural transform between two $\A$-modules is referred to as an $\A$-module morphism. 
We denote by $\ModA$  the category consisting of all $\A$-modules and all $\A$-module morphisms.
Note that  $\ModA$  is obviously an abelian category. 
An $\A$-module  $F$  is called finitely presented if there is an exact sequence $$
\Hom _{\A} (\ \ , X_1) \to \Hom _{\A} (\ \ , X_0) \to F \to 0,  
$$
for some  $X_0, X_1 \in \A$.
We denote by  $\modA$  the full subcategory of  $\ModA$ consisting of all finitely presented $\A$-modules.

\begin{lemma}[Yoneda]
For any  $X \in \A$  and  any $F \in \ModA$, we have the following natural isomorphism:
$$
\Hom _{\ModA} ( \Hom _{\A}(\ \ , X),\  F) \cong F(X).
$$
\end{lemma}

\begin{cor}
An $\A$-module is  projective in  $\modA$  if and only if it is isomorphic to  $\Hom _{\A}(\ \ ,X)$  for some  $X \in \A$. 
\end{cor}

\begin{cor}
The functor  $\A$ to $\modA$  which sends $X$ to $\Hom _{\A}(\ \ ,X)$  is a full embedding.
\end{cor}

Now let  $\C$  be a subcategory of  $\modR$ and let  $\pC$  be the associated stable category. 
Then the category of finitely presented $\C$-modules  $\modC$  and  the category of finitely presented $\pC$-modules  $\modpC$  are defined as in the above course.
Note that for any $F \in \modC$  (resp. $G \in \modpC$)  and for any $X \in \C$ (resp. $\underline{X} \in \pC$),  the abelian group  $F(X)$ (resp. $G(\underline{X})$)  has naturally an $R$-module structure, hence  $F$  (resp. $G$)  is in fact a contravariant additive functor from  $\C$ (resp. $\pC$) to $\modR$.

\begin{rem}
As we stated above there is a natural functor  $\C \to \pC$. 
We can define from this the functor  $\iota : \modpC \to \modC$  by sending 
 $F \in \modpC$ to the composition functor of $\C \to \pC$  with  $F$.
Then it is well known and is easy to prove that  $\iota$  gives an equivalence of categories between  $\modpC$  and the full subcategory of  $\modC$  consisting of all finitely presented $\C$-modules $F$  with  $F(R) = 0$. 
\end{rem}

We prepare the following lemma for a later use.
The proof is straightforward and we leave it to the reader.

\begin{lemma}\label{half-exact}
Let  $0 \to X \to Y \to Z \to 0$  be an exact sequence in  $\modR$. 
Then we have the following.
\begin{itemize}
\item[$(1)$]
The induced sequence  $\pHom_R(W, X) \to  \pHom_R(W, Y) \to \pHom_R(W, Z)$ 
is exact for any $W \in \modR$. 

\item[$(2)$]
If  $\Ext^1_R(Z, R)=0$,   then the induced sequence 
$\pHom_R(Z, W) \to  \pHom_R(Y, W) \to \pHom_R(X, W)$ 
is exact for any $W \in \modR$. 
\end{itemize}
\end{lemma}

\begin{cor}\label{half-exact cor}
Let  $W$  be in $\modR$. 
Then the covariant functor  $\pHom _R(W,  \ \ )$  is a half-exact functor on  $\modR$. 
On the other hand, the contravariant functor  $\pHom_R( \ \ , W)$  is half-exact  on  $\H$.
\end{cor}

\section{Frobenius property of $\modpG$}

\begin{defn}
Let  $\C$   be a subcategory of $\modR$. 

\begin{itemize}
\item[$(1)$]
We say that  $\C$  is closed under kernels of epimorphisms if it satisfies the following condition: 
\begin{itemize}
\item[] 
If $0 \to X \to Y \to Z \to 0$  is an exact sequence in  $\modR$, and if  $Y, Z \in \C$, then  $X \in \C$.
\end{itemize}
(In Quillen's terminology, all epimorphisms from $\modR$ in $\C$ are admissible.) 
\item[$(2)$]
We say that  $\C$  is closed under extension or extension-closed if it satisfies the following condition: 
\begin{itemize}
\item[] 
If $0 \to X \to Y \to Z \to 0$  is an exact sequence in  $\modR$, and if  $X, Z \in \C$, then  $Y \in \C$.
\end{itemize}

\item[$(3)$]
We say that $\C$ is a resolving subcategory if it is extension-closed and closed under kernels of epimorphisms. 

\item[$(4)$] 
We say that  $\C$  is closed under $\Omega$ if it satisfies the following condition: 
\begin{itemize}
\item[] 
If $0 \to X \to P \to Z \to 0$  is an exact sequence in  $\modR$ where  $P$  is a projective module, and if  $Z \in \C$, then  $X \in \C$.
\end{itemize}

Note that for a given $Z \in \C$, the module  $X$  in the above exact sequence is unique up to a projective summand. 
We denote  $X$  by  $\Omega Z$ as an object in $\C$. 
Thus,  $\C$  is closed under $\Omega$  if and only if  $\Omega X \in \C$ whenever  $X \in \C$.

\item[$(5)$]
Similarly to $(3)$, the closedness under $\tr$  is defined. 
Actually, we say that  $\C$  is closed under $\tr$  if  $\tr X \in \pC$ whenever  $X \in \C$. 

\end{itemize}
\end{defn}

Note that the categories $\G$  and  $\H$  satisfy any of the above first four conditions and that  $\G$ is closed under $\tr$. 
We also note the following lemma.  

\begin{lemma}\label{cuke}
Let  $\C$  be a subcategory of  $\modR$.
\begin{itemize}
\item[$(1)$]
If $\C$ is closed under kernels of epimorphisms, then it is closed under $\Omega$.
\item[$(2)$]
If $\C$ is extension-closed and closed under $\Omega$, then it is resolving.
\end{itemize}
\end{lemma}

\begin{pf}
(1) Trivial. 

(2) 
To show that  $\C$  is closed under kernels of epimorphisms, let  $0 \to X \to Y \to Z \to 0$  be an exact sequence in $\modR$ and assume 
that  $Y, Z \in \C$.
Taking a projective cover $P \to Z$ and taking the pull-back diagram, we have the following commutative diagram with exact rows and columns:
$$
\begin{CD}
 @.  @. 0 @. 0  \\
@. @. @VVV @VVV \\
 @.  @. \Omega Z @= \Omega Z  \\
@. @. @VVV @VVV \\
0 @>>> X @>>> E @>>> P @>>> 0 \\
@. @| @VVV @VVV @. \\
0 @>>> X @>>> Y @>>> Z @>>> 0 \\
@. @. @VVV @VVV \\
@.  @. 0 @. 0  \\
\end{CD}
$$
Since  $\C$  is closed under $\Omega$, we have $\Omega Z \in C$.
Then, since  $\C$  is extension-closed, we have  $E \in \C$. 
Noting that the middle row is a split exact sequence, 
we have  $X \in \C$ since  $\C$  is additively closed.
\qed
\end{pf}

We terminologically say that $\cdots \to X_{i+1} \to X_i \to X_{i-1} \to \cdots $ is an exact sequence in a subcategory $\C \subseteq \modR$  if it is an exact sequence in  $\modR$  and such that  $X_i \in \C$  for all $i$.

\begin{prop}\label{resol}
Let  $\C$  be a subcategory of  $\modR$  which is closed under kernels of epimorphisms.
\begin{itemize}
\item[$(1)$]
Then  $\modpC$  is an abelian category with enough projectives.
\item[$(2)$]
For any $F \in \modpC$, there is  a short exact sequence in  $\C$ 
$$
\begin{CD}
0 @>>> X_2  @>>> X_1  @>>> X_0 @>>> 0
 \\
\end{CD}
$$
such that $F$  has a projective resolution of the following type:
$$
\begin{CD}
\cdots 
@>>> \pHom _R ( \ \ , \Omega ^2X_2)|_{\pC}  @>>> \pHom _R ( \ \ , \Omega ^2X_1)|_{\pC}  @>>> \pHom _R ( \ \ , \Omega ^2 X_0)|_{\pC} \\
@>>> \pHom _R ( \ \ , \Omega X_2)|_{\pC}  @>>> \pHom _R ( \ \ , \Omega X_1)|_{\pC}  @>>> \pHom _R ( \ \ , \Omega X_0)|_{\pC} \\
 @>>> \pHom _R ( \ \ , X_2)|_{\pC}  @>>> \pHom _R ( \ \ , X_1)|_{\pC}  @>>> \pHom _R ( \ \ , X_0)|_{\pC}  \\
@>>> F @>>> 0
\end{CD}
$$
\end{itemize}
\end{prop}

\begin{pf}
(1) 
Note that  $\modpC$  is naturally embedded into an abelian category  $\ModpC$.
Let  $\varphi : F \to G$  be a morphism in  $\modpC$. 
It is easy to see just from the definition that $\coker (\varphi) \in \modpC$. 
If we prove that $\ker (\varphi) \in \modpC$, then we see that $\modpC$ is an   abelian category, since it is a full subcategory of the abelian category  $\ModpC$  which is closed under kernels and cokernels.
Now we prove that $\ker (\varphi)$  is finitely presented.

(i) For the first case, we prove it when $F$  and  $G$  are projective.
So  let  $\varphi : \pHom_R( \ \ , X_1) \to \pHom _R(\ \ , X_0)$. 
In this case, by Yoneda's lemma, $\varphi$ is induced from  $\underline{f} : \underline{X_1} \to \underline{X_0}$.
If necessary, adding a projective summand to  $X_1$, we may assume that  $f: X_1 \to X_0$  is an epimorphism in $\modR$. 
Setting  $X_2$  as the kernel of  $f$, we have an exact sequence 
$$
\begin{CD}
0 @>>> X_2 @>>> X_1 @>f>> X_0 @>>> 0. 
\end{CD}
$$
Since  $\C$  is closed under kernels of epimorpshim, we have $X_2 \in \C$. 
Then it follows from Lemma \ref{half-exact} that the sequence 
$$
\begin{CD}
\pHom _R( \ \ , X_2)|_{\pC}  @>{\psi}>> \pHom _R( \ \ , X_1)|_{\pC} @>{\varphi}>> \pHom _R( \ \ , X_0)|_{\pC} 
\end{CD}
$$
 is exact in  $\modpC$. 
Applying the same argument to $\ker (\psi)$, we see that $\ker (\varphi)$  is finitely presented as required. 

(ii) Now we consider a general case.
The morphism $\varphi : F \to G$  induces the following commutative diagram whose horizontal sequences are finite presentations of  $F$  and  $G$: 
$$
\begin{CD}
\pHom _R( \ \ , X_1)|_{\pC} @>{a}>> \pHom _R( \ \ , X_0)|_{\pC} @>{b}>> F @>>> 0 \\
@V{u}VV  @V{v}VV @V{\varphi}VV \\
\pHom _R( \ \ , Y_1)|_{\pC} @>{c}>> \pHom _R( \ \ , Y_0)|_{\pC} @>{d}>> G @>>> 0 \\
\end{CD}
$$
Now we define  $H$  by the following exact sequence:
$$
\begin{CD}
0 @>>> H @>>> \pHom _R( \ \ , X_0)|_{\pC}  \oplus \pHom _R( \ \ , Y_1)|_{\pC} @>{(v, c)}>> \pHom _R( \ \ , Y_0)|_{\pC} 
\end{CD}
$$
It follows from the first step of this proof, we have  $H \in \modpC$.
On the other hand, it is easy to see that there is an exact sequence:
$$
\begin{CD}
\ker (c) \oplus  \pHom _R( \ \ , X_1)|_{\pC}  @>>> H @>>> \ker (\varphi) @>>> 0 \end{CD}
$$
Here, from the first step again, we have  $\ker (c) \in \modpC$.
Now since  $\modpC$  is closed under cokernels in $\ModpC$, we finally have  
$\ker (\varphi)$ as required.

\vspace {6pt}
(2) 
Let  $F$  be an arbitrary object in $\modpC$ with the finite presentation
$$
\begin{CD}
\pHom _R( \ \ , X_1)|_{\pC} @>{\varphi}>> \pHom _R( \ \ , X_0)|_{\pC} @>>> F @>>> 0. \\
\end{CD}
$$
Then, as in the first step of the proof of  (1), we may assume that there is a short exact sequence in  $\C$ 
$$
\begin{CD}
0 @>>> X_2  @>>> X_1  @>f>> X_0 @>>> 0
 \\
\end{CD}
$$
such that  $\varphi$ is induced by  $f$.
Applying Lemma \ref{half-exact}, we have the following exact sequence 
$$
\begin{CD}
\pHom _R( \ \ , X_2)|_{\pC} @>>> \pHom _R( \ \ , X_1)|_{\pC} @>{\varphi}>> \pHom _R( \ \ , X_0)|_{\pC} @>>> F @>>> 0. \\
\end{CD}
$$
Similarly to the proof of Lemma \ref{cuke}(2), taking a projective cover of  $X_0$ and taking the pull-back, we have  the following commutative diagram with exact rows and columns:
$$
\begin{CD}
 @.  @. 0 @. 0  \\
@. @. @VVV @VVV \\
 @.  @. \Omega X_0  @= \Omega X_0  \\
@. @. @VVV @VVV \\
0 @>>> X_2  @>>> E @>>> P @>>> 0 \\
@. @|  @VVV @VVV @. \\
0 @>>> X_2  @>>> X_1  @>>> X_0 @>>> 0 \\
@. @. @VVV @VVV \\
@.  @. 0 @. 0  \\
\end{CD}
$$
 Since the second row is a split exact sequence, we get the exact sequence
$$
\begin{CD}
0 @>>> \Omega X_0 @>>> X_2 \oplus P @>>> X_1 @>>> 0 , 
\end{CD}
$$
where  $P$  is a projective module.
Then it follows from Lemma \ref{half-exact} that there is an exact sequence 
$$
\begin{CD}
\pHom _R( \ \ , \Omega X_0)|_{\pC} @>>>
\pHom _R( \ \ , X_2)|_{\pC} @>>> \pHom _R( \ \ , X_1)|_{\pC}  \\
\end{CD}
$$
Continue this procedure, and we shall obtain the desired projective resolution of  $F$  in $\modpC$.
\qed
\end{pf}

Note that the proof of the first part of Theorem \ref{resol} is completely similar to that of \cite[Lemma (4.17)]{Y}, in which it is proved that $\modpC$  is an abelian category when  $R$  is a Cohen-Macaulay local ring and $\C$  is the category of maximal Cohen-Macaulay modules.

\begin{defn}
A category  $\A$  is said to be a Frobenius category if it satisfies the following conditions:
\begin{itemize}
\item[$(1)$]
$\A$  is an abelian category with enough projectives and with enough injectives.\item[$(2)$] 
All projective objects in $\A$  are injective.
\item[$(3)$] 
All injective objects in $\A$  are projective.
\end{itemize}
Likewise, a category $\A$  is said to be a quasi-Frobenius category if it satisfies the conditions: 
\begin{itemize}
\item[$(1)$]
$\A$  is an abelian category with enough projectives.
\item[$(2)$] 
All projective objects in $\A$  are injective.
\end{itemize}
\end{defn}

\begin{thm}\label{quasi-Frobenius}
Let  $\C$  be a subcategory of  $\modR$  that is closed under kernels of epimorphisms. 
If  $\C \subseteq  \H$  then  $\modpC$  is a quasi-Frobenius category.
\end{thm}

To prove this theorem, we prepare the following lemma.
Here we recall that the full embedding  $\iota : \modpC \to \modC$  is the functor induced by the natural functor  $\C \to \pC$.

\begin{lemma}\label{injective}
Let  $\C$  be a subcategory of  $\modR$  that is closed under kernels of epimorphisms. 
Then the following conditions are equivalent for each  $F \in \modpC$.
\begin{itemize}
\item[$(1)$]
$F$  is an injective object in $\modpC$.
\item[$(2)$]
$\iota F \in \modC$  is half-exact as a functor on $\C$.
\end{itemize}
\end{lemma}

\begin{pf}
As we have shown in the previous proposition, the category  $\modpC$  is an abelian category with enough projectives. 
Therefore an object  $F \in \modpC$  is injective if and only if  $\Ext ^1_{\modpC} (G, F) = 0$  for any  $G \in \modpC$.
But for a given $G \in \modpC$, there is a short exact sequence in $\C$
$$
\begin{CD}
(*) \qquad 0 @>>> X_2  @>>> X_1  @>>> X_0 @>>> 0, 
 \\
\end{CD}
$$
such that  $G$  has a projective resolution 
$$
\begin{CD}
\pHom _R( \ \ , X_2)|_{\pC}  @>>>  \pHom _R( \ \ , X_1)|_{\pC} @>>> \pHom _R( \ \ , X_0)|_{\pC} @>>> G @>>> 0. \\
\end{CD}
$$
Conversely, for any short exact sequence in $\C$ such as  $(*)$, 
the cokernel functor  $G$  of  $\pHom _R (\ \ , X_1)|_{\pC} \to  \pHom _R (\ \ , X_0)|_{\pC}$  is an object of $\modpC$. 
Therefore  $F$  is injective if and only if it satisfies the following condition:

\begin{itemize}
\item[]
The induced sequence 
$$
\begin{CD}
\Hom (\pHom _R (\ \ , X_0)|_{\pC}, F) \to 
\Hom (\pHom _R (\ \ , X_1)|_{\pC}, F) \to 
\Hom (\pHom _R (\ \ , X_2)|_{\pC}, F) 
 \\
\end{CD}
$$
is exact, whenever 
$0  \longrightarrow  X_2   \longrightarrow  X_1  \longrightarrow  X_0  \longrightarrow  0$ is a short exact sequence in $\C$.

\end{itemize}

It follows from Yoneda's lemma that this is equivalent to saying that 
$F(X_0)  \longrightarrow  F(X_1)   \longrightarrow  F(X_2)$ 
is exact whenever 
$0  \longrightarrow  X_2   \longrightarrow  X_1  \longrightarrow  X_0  \longrightarrow  0$ 
is a short exact sequence in $\C$.
This exactly means that  $\iota F$ is half-exact as a functor on  $\C$. 
\qed
\end{pf}

\begin{pf}
Now we proceed to the proof of Theorem \ref{quasi-Frobenius}. 
We have already shown that  $\modpC$ is an abelian category with enough projectives.
It remains to show that any projective module $\pHom _R ( \ \ , X)|_{\pC}$ ($X \in \C$)  is an injective object in $\modpC$.
Since  $\C$  is a subcategory of $\H$,  it follows from Corollary \ref{half-exact cor} that  $\pHom_R(\ \ , X)|_{\C} = \iota ( \pHom _R(\ \ , X) |_{\pC})$  is a half-exact functor, hence it is injective by the previous lemma.
\qed
\end{pf}

Before stating the next theorem, we should remark that the syzygy functor  $\Omega$  gives an automorphism on $\pG$.

\begin{thm}\label{Frob}
Let  $\C$  be a subcategory of $\modR$. 
And suppose the following conditions. 
\begin{itemize}
\item[$(1)$]
$\C$  is a resolving subcategory of $\modR$. 
\item[$(2)$]
$\C \subseteq \H$. 
\item[$(3)$]
The functor  $\Omega : \pC \to \pC$ yields  a surjective map on the set of isomorphism classes of the objects in $\pC$. 
\end{itemize}
Then  $\modpC$  is a Frobenius category. 
In particular,  $\modpG$ is a Frobenius category.
\end{thm}

\begin{pf} 
Since  $\C$  is subcategory of  $\H$  that is closed under kernels of epimorphisms in  $\modR$,  $\modpC$ is a quasi-Frobenius category by the previous theorem. It remains to prove that  $\modpC$  has enough injectives and all injectives are projective.

(i) 
For the first step of the proof we show that each $\pC$-module  $F \in \modpC$  can be embedded into a projective $\pC$-module  $\pHom _R( \ \ , Y)|_{\pC}$  for some $Y \in \C$.

In fact, as we have shown in Proposition \ref{resol},  for a given  $F \in \modpC$, there is a short exact sequence in $\C$ 
$$
\begin{CD}
0 @>>> X_2 @>>> X_1 @>>> X_0 @>>> 0
\end{CD}
$$
such that  $F$  has a projective resolution 
$$
\begin{CD}
\pHom _R( \ \ , X_1)|_{\pC} @>>> \pHom _R( \ \ , X_0)|_{\pC} @>>> F @>>> 0. \\
\end{CD}
$$
From the assumption (3), there is an exact sequence in  $\C$ 
$$
\begin{CD}
0 @>>> X_2 @>>> P  @>>> Y @>>> 0, 
\end{CD}
$$
where  $P$  is a projective module. 
Then, similarly to the argument in the proof of Lemma \ref{cuke}, just taking the push-out, we have a commutative diagram with exact rows and exact columns :
$$
\begin{CD}
@. 0 @. 0 @.  @. \\
@. @VVV @VVV @. @. \\
0 @>>> X_2 @>>> X_1 @>>> X_0 @>>> 0 \\
@. @VVV @VVV @| @. \\
0 @>>> P  @>>> E  @>>> X_0 @>>> 0 \\
@. @VVV @VVV @.  @. \\
 @.  Y  @=  Y  @.  \\
@. @VVV @VVV @.  @. \\
@. 0 @. 0 @. \\
\end{CD}
$$
Since  $\Ext ^1 _R (X_0, P) =0$, we should note that the second row is splittable. 
Hence we have a short exact sequence of the type 
$$
\begin{CD}
0 @>>> X_1 @>>> X_0 \oplus P  @>>> Y @>>> 0. 
\end{CD}
$$
Therefore we have from Lemma \ref{half-exact} that there is an exact sequence 
$$
\begin{CD}
\pHom _R(\ \ , X_1)|_{\pC}  @>>> \pHom _R(\ \ , X_0)|_{\pC}  @>>> \pHom _R(\ \ ,  Y)|_{\pC}.  
\end{CD}
$$
Thus  $F$  can be embedded into  $\pHom _R(\ \ , Y)|_{\pC}$  as desired. 
\par
(ii) Since  all projective modules in  $\modpC$  are injective, it follows from  (i) that  $\modpC$  has enough injectives.

\par
(iii) To show that every injective module in $\modpC$ is projective, let  $F$  be an injective  $\pC$-module in $\modpC$.
 By (i), $F$  is a $\pC$-submodule of  $\pHom _R(\ \ , Y)|_{\pC}$ for some  $Y \in \C$, hence  $F$  is a direct summand of  $\pHom _R(\ \ , Y)|_{\pC}$. 
Since  it is a summand of a projective module,  $F$  is projective as well. 
\qed
\end{pf}

\section{Characterizing subcatgories of $\H$}

In this section we always assume that $R$  is a henselian local ring with maximal ideal  $\m$ and with the residue class field  $k =R/ \m$. 
In the following, what we shall need from this assumption is the fact that 
  $X \in \modR$  is indecomposable only if  $\End _R(X)$  is a (noncommutative)  local ring.  
In fact we can show the following lemma.

\begin{lemma}\label{KS thm}
Let  $(R, \m)$ be a henselian local ring. 
Then  $\modpC$  is a Krull-Schmidt category  
for any subcategory  $\C \subseteq \modR$.
\end{lemma}

\begin{pf}
We have only to prove that  $\End _{\modpC} (F)$  is a local ring for any indecomposable $\pC$-module $F \in \modpC$. 

First we note that $\End _{\modpC} (F)$ is a module-finite algebra over  $R$. 
In fact, since there is a finite presentation 
$$
\begin{CD}
\pHom _R (\ \ , X_1)|_{\pC} @>>> \pHom _R (\ \ , X_0)|_{\pC}  @>>> F  @>>> 0,  \\
\end{CD}
$$ 
$F(X_0)$  is a finite $R$-module. 
On the other hand, taking the dual by  $F$ of the above sequence and using Yoneda's lemma, we can see that there is an exact sequence of $R$-modules 
$$
\begin{CD}
0 @>>> \End _{\modpC} (F) @>>> F(X_0) @>>> F(X_1).  \\
\end{CD}
$$
As a submodule of a finite module, $\End _{\modpC} (F)$  is finite over  $R$. 

Now suppose that  $\End _{\modpC} (F)$  is not a local ring. 
Then there is an element  $e \in \End _{\modpC} (F)$  such that both $e$ and  $1-e$ are nonunits. 
Let  $\overline{R}$  be the image of the natural ring homomorphism  $R \to \End _{\modpC} (F)$. 
 We consider the subalgebra  $\overline{R}[e]$  of $\End _{\modpC} (F)$. 
 Since  $\overline{R}[e]$ is a commutative $R$-algebra which is finite over  $R$, it is also henselian. 
 Since  $e, 1-e \in \overline{R}[e]$ are both nonunits, 
  $\overline{R}[e]$ should be decomposed into a direct product of rings, in particular it contains a nontrivial idempotent. 
 This implies that  $\End _{\modpC} (F)$  contains a nontrivial idempotent, hence that  $F$  is decomposable in  $\modpC$. 
\qed
\end{pf}

We claim the validity of the converse of Theorem \ref{quasi-Frobenius}. 

\begin{thm}\label{characterize H}
Let   $\C$   be a resolving subcategory of $\modR$.
Suppose that $\modpC$  is a quasi-Frobenius category. 
Then  $C \subseteq \H$. 
\end{thm}

In a sense  $\H$  is the largest resolving subcategory $\C$  of $\modR$  for which  $\modpC$ is a quasi-Frobenius category.

\begin{pf}
Let  $X$  be an indecomposable nonfree module in  $\C$. 
It is sufficient to prove that  $\Ext^1_R (X, R) =0$. 

In fact, if it is true for any  $X \in \C$, then  $\Ext ^1 _R (\Omega ^iX, R) =0$ for any $i \geq 0$, since  $\C$  is closed under $\Omega$.
This implies that  $\Ext ^{i+1} _R (X, R) =0$  for $i \geq 0$, hence  $X \in \H$.

Now we assume that $\Ext ^1_R(X, R) \not= 0$  to prove the theorem by contradiction. 
Let  $\sigma$  be a nonzero element of  $\Ext ^1_R(X, R)$  that corresponds to the nonsplit extension 
$$
\begin{CD}
\sigma \ : \  0 @>>> R @>>> Y @>p>> X @>>> 0.
\end{CD}
$$
Since  $\C$  is closed under extension, we should note that  $Y \in \C$.
On the other hand, noting that  $R$ is the zero object in $\pC$, we have from Lemma \ref{half-exact} that there is an exact sequence of  $\pC$-modules 
$$
\begin{CD}
0  @>>> \pHom _R(\ \ , Y)|_{\pC} @>{p_*}>>\pHom _R(\ \ , X)|_{\pC}.
\end{CD}
$$
Since  $\modpC$  is a quasi-Frobenius category, this monomorphism is a split one, hence  $\pHom _R(\ \ , Y)|_{\pC}$  is a direct summand of  $\pHom _R(\ \ , X)|_{\pC}$  through  $p_*$. 
Since  the embedding  $\pC \to \modpC$  is full, and since we assumed $X$  is indecomposable, $\pHom _R(\ \ , X)|_{\pC}$  is indecomposable in $\modpC$ as well. 
Hence  $p_*$  is either an isomorphism or $p_* =0$.

We first consider the case that  $p_*$  is an isomorphism. 
In this case, we can take a morphism  $q \in \Hom _R(X, Y)$  such that  $q_* : \pHom _R(\ \ , X)|_{\pC} \to \pHom _R(\ \ , Y)|_{\pC}$  is the inverse of  $p_*$.  Then  $p_*q_* =(pq)_*$  is the identity on  $\pHom _R(\ \ , X)|_{\pC}$. 
Since $\pEnd_R (X) \cong \End (\pHom _R(\ \ , X)|_{\pC})$, we see that 
$\underline{pq} = \underline{1}$ in $\pEnd_R(X)$. 
Since  $\End _R (X)$  is a local ring and since  $\pEnd_R (X)$  is a residue ring of  $\End _R(X)$, we see that  $pq \in \End _R(X)$  is a unit. 
This shows that the extension $\sigma$ splits, 
and this is a contradiction. 
Hence this case never occurs. 

 As a result we have that $p_* = 0$, which implies  $\pHom _R (\ \ , Y)|_{\pC}= 0$.
Since  the embedding  $\pC \to \modpC$  is full, this is equivalent to saying that  $Y$  is a projective, hence free, module. 
Thus it follows from the extension $\sigma$ that $X$  has projective dimension exactly one. 
(We have assumed that  $X$  is nonfree.)
Thus the extension $\sigma$ should be a minimal free resolution of  $X$: 
$$
\begin{CD}
0 @>>> R @>{\alpha}>> R^r @>p>> X @>>> 0.
\end{CD}
$$
where  $\alpha = (a_1, \cdots , a_r)$  is a matrix with entries in $\m$. 
Now let  $x \in \m$  be any element and let us consider the extension corresponding to  $x\sigma \in \Ext ^1_R (X, R)$. 
By making push-out, we obtained this extension as the second row in the following commutative diagram with exact rows : 
$$
\begin{CD}
\sigma \ @.: \ 0 @>>> R @>{\alpha}>> R^r @>p>> X @>>> 0 \\
@.  @.  @VxVV  @VVV @| \\ 
x\sigma \ @.: \ 0 @>>> R @>>> Z @>{p'}>> X @>>> 0 
\end{CD}
$$
Note that there is an exact sequence 
$$
\begin{CD}
 R  @>{(x , \alpha)}>> R \oplus R^r @>>> Z @>>> 0, 
\end{CD}
$$
where all entries of the matrix $(x, \alpha)$ are in $\m$. 
Thus the $R$-module  $Z$ is not free, and thus $p'_* : \pHom_R(\ \ , Z)|_{\pC} \to \pHom_R(\ \ , X)|_{\pC}$ is a nontrivial monomorphism.  
Then, repeating the argument in the first case to the extension  $x\sigma$, 
we must have   $x \sigma = 0$  in  $\Ext^1_R(X,R)$. 
Since this is true for any  $x \in \m$  and for any $\sigma \in \Ext^1_R(X, R)$, we obtain that  $\m \Ext ^1_R(X, R) = 0$. 
On the other hand, by computation, we have  $\Ext^1_R (X, R) \cong R/(a_1, \ldots , a_r)$, and hence we must have  $\m = (a_1, \ldots , a_r)R$. 
Since the residue field  $k$  has a free resolution of the form 
$$
\begin{CD}
 R^r   @>{{}^t(a_1, \ldots , a_r) }>> R @>>> k @>>> 0, 
\end{CD}
$$
comparing this with the extension $\sigma$, we have that $X \cong \tr k$. 
What we have proved so far is the following :

\begin{itemize}
\item[]
 Suppose there is an indecomposable nonfree module  $X$  in $\C$ which satisfies  $\Ext^1_R (X, R) \not= 0$. 
Then  $X$ is isomorphic  to $\tr k$  as an object in  $\pC$  and  $X$  has projective dimension one. 
\end{itemize}

If  $R$  is a field then the theorem is obviously true. 
So we assume that the local ring  $R$  is not a field. 
 Then we can find an indecomposable $R$-module $L$  of length $2$ and a nonsplit exact sequence 
$$
\begin{CD}
0 @>>>  k  @>>> L @>>> k @>>> 0. 
\end{CD}
$$
 Note that  $L = R/I$  for some $\m$-primary ideal $I$. 
 Note also that  $\depth R \geq 1$, since there is a module of projective dimension one. 
 Therefore there is no nontrivial $R$-homomorphism from  $L$  to  $R$, and thus  we have  $\pEnd _R (L) = \End _R(L) \cong R/I$. 
 Similarly we have  $\pEnd _R (k) \cong k$. 

It also follows from  $\Hom _R (k , R)= 0$  that there is an exact sequence of the following type: 
$$
\begin{CD}(*) \quad 
0 @>>>  \tr k  @>>> \tr L \oplus P @>>> \tr k @>>> 0, 
\end{CD}
$$
where $P$ is a suitable free module. 
 Since $X = \tr k$  is in $\C$, and since  $\C$  is extension-closed, we have 
 $\tr L \in \C$ as well. 
 Note that  $\tr$  is a duality on  $\pmodR$, hence we see that  $\pEnd _R (\tr L ) \cong \pEnd _R (L)$  that is a local ring. 
As a consequence we see that  $\tr L$  is indecomposable in  $\pC$. 
 
We claim that  $\tr L$  has projective dimension exactly one, hence in particular $\Ext^1_R (\tr L, R) \not= 0$. 
In fact, we have from (*) that $\tr L$  has projective dimension at most one as well as $X = \tr k$. 
If $\tr L$  were free then  $X = \tr k$  would be its own first syzygy from (*) and hence free because $X$ has projective dimension one. 
But this is a contradiction. 

Thus it follows from the above claim that  $\tr L$  is isomorphic to $\tr k$  in $\pmodR$. 
 Taking the transpose again, we finally have  that $L$  is isomorphic to  $k$  in $\pmodR$. 
 But this is absurd, because  $\pEnd _R (L) \cong R/I$  and  $\pEnd _R (k) \cong k$. 
Thus the proof is complete. 
\qed
\end{pf}

\begin{cor}
The following conditions are equivalent for a henselian local ring  $R$.
\begin{itemize}
\item[$(1)$]
$\mod(\pmodR)$ is a quasi-Frobenius category.
\item[$(2)$]
$\mod(\pmodR)$ is a Frobenius category.
\item[$(3)$]
$R$  is an artinian Gorenstein ring.
\end{itemize}
\end{cor}

\begin{pf}
$(3) \Rightarrow (2)$: 
If $R$  is an artinian Gorenstein ring, then we have  $\modR = \G$.
Hence this implication follows from Theorem \ref{Frob}.

$(2) \Rightarrow (1)$: Obvious.

$(1) \Rightarrow (3)$: 
Suppose that $\mod(\pmodR)$ is a quasi-Frobenius category. 
Then, by Theorem \ref{characterize H}, any indecomposable $R$-module $X$ is in  $\H$. 
 In particular, the residue field  $k = R/\m$  is in  $\H$, hence by definition,  $\Ext ^i_R (k, R) = 0$  for any $i > 0$. 
 This happens only if  $R$  is an artinian Gorenstein ring.   
\qed
\end{pf}

\section{Characterizing subcatgories of $\G$}

\begin{lemma}\label{main lemma}
Let  $R$  be a henselian local ring and let  $\C$ be an extension-closed  subcategory of  $\modR$. 
For objects  $\underline{X}, \underline{Y} \in \pC$, we assume the following:
\begin{itemize}
\item[$(1)$]
There is a monomorphism $\varphi$  in $\ModpC$:
$$
\varphi : \pHom _R( \ \ , Y)|_{\pC}  \to \Ext^1( \ \ , X)|_{\pC}
$$
\item[$(2)$]
$\underline{X}$  is indecomposable in $\pC$.
\item[$(3)$]
$\underline{Y} \not\cong 0$  in $\pC$.
\end{itemize}
Then the module $X$  is isomorphic to a direct summand of  $\Omega Y$. 
\end{lemma}

\begin{pf}
We note from Yoneda's lemma that there is an element  $\sigma \in \Ext ^1_R (Y, X)$  which corresponds to the short exact sequence :
$$
\begin{CD}
\sigma : 0 @>>> X @>a>> L @>p>> Y @>>> 0 \\
\end{CD}
$$   
such that  $\varphi$  is induced by $\sigma$  as follows: 

For any  $W \in \C$ and for any  $\underline{f} \in \pHom _R (W, Y)$, 
 consider the pull-back diagram to get the following commutative diagram with exact rows: 
$$
\begin{CD}
0 @>>> X @>>> L @>>> Y @>>> 0 \\
@.  @|   @AAA  @AfAA \\
0 @>>> X @>>> E @>>> W @>>> 0 \\ 
\end{CD}
$$   
Setting the second exact sequence as  $f^* \sigma$, we have  $\varphi (\underline{f}) = f^* \sigma$. 

Note that  $L \in \C$, since  $\C$  is extension-closed. 
Also note that  $\sigma$ is nonsplit. 
In fact, if it splits, then  $\varphi$  is the zero map, hence  $\pHom_R(\ \ ,Y)|_{\pC} = 0$ from the assumption. 
Since the embedding  $\pC \to \modpC$  is full, this implies  that  $\underline{Y} = \underline{0}$  in $\pC$, which is a contradiction. 

Now let  $P \to Y$  be a surjective $R$-module homomorphism where  $P$  is a projective module. 
 Then there is a commutative diagram with exact rows: 
$$
\begin{CD}
 0 @>>> X @>a>> L @>p>> Y @>>> 0 \\
@.  @AgAA   @AhAA  @| \\
 0 @>>> \Omega Y @>{\alpha}>> P @>{\pi}>> Y @>>> 0 \\ 
\end{CD}
$$   
Note that the extension $\sigma$  induces the exact sequence of  $\pC$-modules:
$$
\begin{CD}
\pHom _R (\ \ , L)|_{\pC}  @>{p_*}>> \pHom _R (\ \ , Y)|_{\pC}  @>{\varphi}>> \Ext ^1 _R (\ \ , X)|_{\pC} 
\end{CD}
$$   
 Since $\varphi$ is a monomorphism, the morphism 
$\pHom _R (\ \ , L)|_{\pC}  \to  \pHom _R (\ \ , Y)|_{\pC}$ 
is the zero morphism. 
In particular,  the map $\underline{p} \in \pHom _R (L, Y)$ is the zero element by Yoneda's lemma. 
(Note that we use the fact  $L \in \C$ here.) 
This is equivalent to saying that  $p : L \to Y$ factors through a projective module, hence that it  factors through the map $\pi$. 
As a consequence, there are  maps  $k : L \to P$  and  $\ell : X \to \Omega Y$  which make the following diagram commutative: 
$$
\begin{CD}
 0 @>>> X @>a>> L @>p>> Y @>>> 0 \\
@.  @V{\ell}VV   @V{k}VV   @|  \\
 0 @>>> \Omega Y @>{\alpha}>> P @>{\pi}>> Y @>>> 0 \\ 
\end{CD}
$$   
Then, since  $p(1-hk)=0$, there is a map $b : L  \to X$ such that  $1 -hk =ab$.  Likewise, since  $\pi (1-kh) =0$, there is a map $\beta : P \to \Omega Y$ such that  $1 -kh = \alpha \beta$.
Note that there are  equalities: 
$$
a(1-g \ell) = a -ag \ell = a - h \alpha \ell = a - hka = (1-hk)a = aba
$$
Since  $a$  is a monomorphism, we hence have $1 - g \ell =  ba$. 
 Thus we finally obtain the equality  $1 = ba + g \ell$  in the local ring   $\End _R (X)$. 

Since  $\sigma$  is a nonsplit sequence, $ba \in \End _R(X)$  never be a unit, 
 and it follows that  $g \ell$  is a unit in  $\End _R (X)$. 
 This means that the map  $g : \Omega Y \to X$  is a split epimorphism, hence  $X$  is isomorphic to a direct summand of  $\Omega Y$ as desired. 
\qed
\end{pf}

\begin{thm}\label{moduloAR}
Let  $R$  be a henselian local ring.
Suppose that 
\begin{itemize}
\item[$(1)$]
$\C$  is a resolving subcategory of  $\modR$.
\item[$(2)$]
$\modpC$ is a Frobenius category.
\item[$(3)$]
There is no nonprojective module  $X \in \C$  with  $\Ext^1_R(\ \ ,X)|_{\pC} = 0$.
\end{itemize}
Then  $\C \subseteq \G$.
\end{thm}

\begin{pf}
As the first step of the proof, we prove the following:
\begin{itemize}
\item[(i)]
 For a nontrivial indecomposable object  $\underline{X} \in \pC$, there is an object  $\underline{Y} \in \pC$  such that  $X$  is isomorphic to a direct summand of  $\Omega Y$.
\end{itemize}

To prove this, let  $\underline{X} \in \pC$  be nontrivial and indecomposable.
Consider the $\pC$-module  $F : = \Ext^1_R (\ \ , X)|_{\pC}$. 
 The third assumption assures us that $F$  is a nontrivial  $\pC$-module.
 Hence there is an indecomposable module  $W \in \pC$ such that 
 $F(W) \not =0$. 
Take a nonzero element $\sigma$ in  $F(W) = \Ext^1_R (W, X)$ that corresponds to an exact sequence
$$
\begin{CD}
0 @>>> X @>>> E @>>> W @>>> 0. \\
\end{CD}
$$
Note that  $E \in \C$,  since  $\C$  is extension-closed.
Then we have an exact sequence of  $\pC$-modules 
$$
\begin{CD}
\pHom _R (\ \ , E)|_{\pC} @>>> \pHom _R (\ \ , W)|_{\pC} @>{\varphi}>> \Ext ^1 _R (\ \ , X)|_{\pC} \\
\end{CD}
$$
We denote by  $F_{\sigma}$  the image of  $\varphi$.  
 Of course,  $F _{\sigma}$  is a nontrivial $\pC$-submodule of  $F$   which is finitely presented.
 Since we assume that  $\modpC$  is a Frobenius category, we can take a minimal  injective hull of  $F_{\sigma}$  that is projective as well, i.e. 
 there is a monomorphism  $i: F _{\sigma} \to  \pHom_R (\ \ , Y)|_{\pC}$  for some  $Y \in \pC$ that is an essential extension.  
 Since  $\Ext _R^1(\ \ , X)|_{\pC}$  is half-exact as a functor on  $\C$, 
 we can see by a similar method to that in the proof of Lemma \ref{injective} that  $\Hom (\ \ , \Ext ^1_R(\ \ ,X))$  is an exact functor on  $\modpC$. 
 It follows from this that  the natural embedding  $F_{\sigma} \to F$  can be enlarged to the morphism  $g : \pHom_R (\ \ , Y)|_{\pC} \to F$. 
Hence there is a commutative diagram 
$$
\begin{CD}
F_{\sigma}  @>{\subset}>{i}> \pHom _R(\ \ , Y)|_{\pC} \\
@V{\cap}VV    @VgVV  \\
F   @= \Ext ^1 _R(\ \ , X). \\
\end{CD}
$$ 
Since  $i$ is an essential extension, we see that  $\ker \ g =0$, hence we have   $\pHom _R (\ \ , Y)|_{\pC}$  is a submodule of  $F$. 
 Hence by the previous lemma we see that  $X$  is isomorphic to a direct summand of  $\Omega Y$. 
 Thus the claim (i)  is proved.

Now we prove the theorem. 
Since  $\modpC$  is a quasi-Frobenius category, we know from Theorem 
 \ref{characterize H} that  $\C \subseteq \H$. 
 To show  $\C \subseteq \G$,  let  $\underline{X}$  be a nontrivial indecomposable module in $\pC$.  
We want to prove that  $\Ext ^i_R (\tr X, R) = 0$ for $i >0$. 
It follows from the claim (i) that there is  $Y \in \pC$  such that 
 $X$  is a direct summand of  $Y$. 
Note that $Y \in \H$. 
From the obvious sequence 
$$
\begin{CD}
0 @>>> \Omega Y @>>> P @>>> Y @>>> 0 \\
\end{CD}
$$ 
with  $P$ being a projective module, it is easy to see that 
there is an exact sequence of the type 
$$
\begin{CD}
(*) \qquad 0 @>>> \tr Y @>>> P' @>>> \tr \Omega Y @>>> 0,  \\
\end{CD}
$$ 
where  $P'$  is projective.
Since  $\Omega Y$  is a torsion-free module, it is obvious that 
$\Ext^1_R (\tr \Omega Y, R)=0$. 
 Since  $X$  is a direct summand of  $\Omega Y$, we have that 
$\Ext ^1_R (\tr X, R) = 0$ as well.
  This is true for any indecomposable module in  $\pC$, hence for each indecomposable summand of  $Y$. 
Therefore we have  $\Ext^1_R(\tr Y, R) =0$. 
  Then it follows from this together with  (*)  that 
 $\Ext ^2 _R(\tr \Omega Y, R) =0$, hence we have that  $\Ext ^2 _R(\tr X, R) =0$. 
 Continuing this procedure, we can show  $\Ext_R^i (\tr X, R) =0$  for any $i >0$.
\qed
\end{pf}

\begin{rem}
We conjecture that  $\G$  should be the largest resolving subcategory  $\C$  of  $\modR$  such that  $\modpC$  is a Frobenius category.

Theorem \ref{moduloAR}  together with Theorem \ref{Frob} say that this is true modulo Auslander-Reiten conjecture:

\begin{itemize}
\item[$(AR)$] 
If  $\Ext^i_R(X, X \oplus R) = 0$  for any $i>0$ then  $X$ should be projective.\end{itemize}

In fact, if the conjecture (AR)  is true, then the third condition of the previous theorem is automatically satisfied. 
\end{rem}

\begin{defn}
Let  $\A$  be any additive category. 
We denote by  $\ind (\A)$  the set of nonisomorphic modules which represent  all the isomorphism classes of indecomposable objects in  $\A$. 
 If  $\ind (\A)$  is a finite set, then we say that  $\A$  is a category of finite type. 
\end{defn}

The following theorem is a main theorem of this paper, which claims that any resolving subcategory of finite type in  $\H$  are contained in  $\G$.

\begin{thm}\label{characterize G}
Let  $R$  be a henselian local ring and let  $\C$  be a subcategory of $\modR$ which satisfies the following conditions.
\begin{itemize}
\item[$(1)$]
$\C$ is a resolving subcategory of $\modR$. 
\item[$(2)$]
$\C \subseteq \H$.
\item[$(3)$]
$\C$  is of finite type.
\end{itemize}
Then, $\modpC$  is a Frobenius category and  $\C \subseteq \G$.
\end{thm}

\begin{lemma}\label{finiteness}
Let  $\C$  be a subcategory of  $\modR$  and suppose that $\C$  is of finite representation type.  
Then the following conditions are equivalent for a contravariant additive functor   $F$  from  $\pC$  to  $\modR$. 
\begin{itemize}
\item[$(1)$]  $F$  is finitely presented, i.e. $F \in \modpC$. 
\item[$(2)$]  $F(W)$  is a finitely generated $R$-module for each  $W \in \ind (\pC)$. 
\end{itemize}
\end{lemma}

\begin{pf}
The implication  $(1) \Rightarrow (2)$  is trivial from the definition.
 We prove  $(2) \Rightarrow (1)$. 
 For this, note that for each  $X, W \in \ind (\pC)$  and for each $\underline{f} \in \pHom _R (X, W)$, the induced map  $F(f) : F(W) \to F(X)$  is an $R$-module homomorphism and  satisfies that  $F(af)=aF(f)$  for  $a \in R$. 
Therefore the $\pC$-module homomorphism 
$$
\varphi_{W} : \pHom _R(\ \ , W) \otimes _R F(W)  \to F
$$
which sends $f\otimes x$  to  $F(f)(x)$  is well-defined.
 Now let  $\{ W_1 , \cdots , W_m \}$  be the complete list of elements in  $\ind (\pC)$.  
 Then the $\pC$-module homomorphism 
$$
\Phi = \oplus _{i=1}^m \varphi _{W_i} : \oplus _{i=1}^m \pHom _R(\ \ , W_i) \otimes _R F(W_i)  \to F 
$$
is defined, and it is clear that $\Phi$  is an epimorphism in  $\ModpC$. 
Therefore  $F$  is finitely generated, and this is true for $\ker (\Phi)$  as well. 
Hence  $F$  is finitely presented. 
\qed
\end{pf}

\begin{lemma}\label{AR}
Let  $\C$  be a subcategory of  $\modR$  that is of finite type.  
And let  $\underline{X}$  be a nontrivial indecomposable module in  $\pC$. 
 Suppose that  $\C$  is closed under kernels of epimorphisms. 
 Then $\C$  admits an AR-sequence ending in  $X$,  that is, 
there is a nonsplit exact sequence in $\C$ 
$$
\begin{CD}
0 @>>> \tau X @>>> L  @>p>> X @>>> 0,  \\
\end{CD}
$$ 
such that  for any indecomposable  $Y \in \C$  and for any morphism  $f : Y \to  X$  which is not a split epimorphism, there is a morphism  $g : Y \to L$ that makes the following diagram commutative. 
$$
\begin{CD}
L  @>p>> X   \\
@AgAA @AfAA \\
Y @= Y \\
\end{CD}
$$
\end{lemma}

\begin{pf}
Let  $\rad \ \pHom _R ( \ \ , X)|_{\pC} $  be the radical functor of  $\pHom _R( \ \ , X)|_{\pC}$, i.e. for each  $\underline{W} \in \ind (\pC)$,  if  $\underline{W} \not\cong \underline{X}$ then  $\rad \ \pHom_R(W, X) = \pHom _R (W, X)$, on the other hand, if  $\underline{W} = \underline{X}$ then  $\rad \ \pHom _R (X,X)$  is the unique maximal ideal of  $\pEnd _R(X)$. 
 Since  $\rad \ \pHom _R  ( \ \ , X)|_{\pC} $ is a $\pC$-submodule of  $\pHom _R ( \ \ , X) |_{\pC}$, it follows from the previous lemma that 
 $\rad \ \pHom _R  ( \ \ , X)|_{\pC} $  is finitely presented, hence 
there is an $\underline{L} \in \pC$  and a morphism  $\underline{p} : \underline{L} \to \underline{X}$  such that  $p_* : \pHom _R  ( \ \ , L )|_{\pC} \to \rad \ \pHom _R  ( \ \ , X)|_{\pC}$  is an epimorphism. 
Adding a projective summand to  $L$ if necessary, we may assume that the $R$-module homomorphism  $p : L \to X$  is surjective. 
Setting  $\tau X = \ker (p)$, we see that  $\tau X \in \C$, since  $\C$  is closed under kernels of epimorphisms. 
 And it is clear that the obtained sequence  $0 \to \tau X \to L \to X \to 0$  satisfies the required condition to be an AR-sequence. 
\qed
\end{pf}

See \cite{A} and \cite{AR}  for the detail of AR-sequences. 

\begin{lemma}\label{finite length}
Let  $\C$  be a resolving subcategory of  $\modR$. 
 Suppose that  $\C$  is  of finite type. 
 Then for any  $X, Y \in  \C$,  the $R$-module 
$\Ext^1_R (X, Y)$  is of finite length. 
\end{lemma}

\begin{pf}
It is sufficient to prove the lemma in the case that  $X$  and  $Y$  are indecomposable. 
 For any  $x \in \m$  and for any $\sigma \in \Ext ^1_R (X, Y)$, it is enough to show that  $x^n \sigma = 0$  for a large integer $n$. 
 
Now suppose that  $x^n \sigma \not= 0$  for any integer $n$, and we shall show a contradiction.
 Let us take an AR-sequence ending in $X$ as in the previous lemma 
$$
\begin{CD}
\alpha : \ 0 @>>> \tau X  @>>> L  @>p>> X @>>> 0,  \\
\end{CD}
$$ 
and  a short exact sequence that corresponds to each  $x^n \sigma \in \Ext ^1 _R (X, Y)$
$$
\begin{CD}
x^n \sigma :  \ 0 @>>> Y  @>>> L_n  @>p_n>> X @>>> 0.  \\
\end{CD}
$$ 
Since  $p_n$ is not a split epimorphism, the following commutative diagram is induced:
$$
\begin{CD}
0 @>>> Y  @>>> L_n  @>p_n>> X @>>> 0  \\
@. @V{h_n}VV  @VVV @| \\ 
0 @>>> \tau X @>>> L @>p>> X @>>> 0 \\
\end{CD}
$$ 
The morphism  $h_n$ induces an $R$-module map 
$$
\begin{CD}
(h_n)_* : \Ext^1_R (X, Y) @>>> \Ext ^1_R (X, \tau X) 
\end{CD}
$$ 
which sends  $x^n\sigma$  to the AR-sequence  $\alpha$. 
Since  $(h_n)_*$  is $R$-linear, we have 
$\alpha = x^n (h_n)_*(\sigma) \in x^n \Ext ^1_R(X, \tau X)$. 
 Note that this is true for any integer $n$  and that 
  $\cap _{i=1} ^{\infty} x^n \Ext ^1_R(X, \tau X) = (0)$. 
 Therefore we must have  $\alpha = 0$. 
This contradicts to that $\alpha$  is a nonsplit exact sequence. 
\qed
\end{pf}

\begin{rem}
Compare the proof of the above lemma with that in \cite[Theorem (3.4)]{Y}.
\end{rem}

\vspace{12pt}
 
Now we proceed to the proof of Theorem \ref{characterize G}. 
For this, let  $\C$  be a subcategory of $\modR$ that satisfies three conditions as in the theorem. 
The proof will be done step by step. 

For the first step we show that 

\begin{itemize}
\item[(Step 1)] the category  $\modpC$  is a quasi-Frobenius category. 
\end{itemize}

This has been proved in Theorem \ref{quasi-Frobenius}, since  $\C$  is a resolving subcategory of  $\H$. 
\qed

Now we prove the following.  

\begin{itemize}
\item[(Step 2)] Any $\pC$-module  $F \in \modpC$  can be embedded in an injective $\pC$-module of the form  $\Ext^1_R (\ \ , X)|_{\pC}$  for some  $\underline{X} \in \pC$. 
In particular,  $\modpC$  has enough injectives. 
\end{itemize}

\begin{pf}
As we have shown in Lemma \ref{resol} that for a given  $F \in \modpC$, there is a short exact sequence  in  $\pC$ 
$$
\begin{CD}
 0 @>>> X_2  @>>> X_1  @>>> X_0  @>>> 0  \\
\end{CD}
$$ 
such that a projective resolution of  $F$  in  $\modpC$  is given as in Lemma \ref{resol}(2). 
It is easy to see from the above exact sequence that there is an exact sequence $$
\begin{CD}
\cdots @>>> \pHom_R(\ \ , X_1)|_{\pC}  @>>>  \pHom_R(\ \ , X_0)|_{\pC}  @>>> 
 \Ext ^1_R(\ \ , X_2)|_{\pC} @>>> \cdots 
\end{CD}
$$ 
 Hence there is a monomorphism  $F  \to   \Ext ^1_R(\ \ , X_2)|_{\pC}$. 
 Note that  $\Ext ^1_R(W, X_2)$  is a finitely generated $R$-module for each  $\underline{W} \in \ind (\pC)$. 
Hence it follows from Lemma \ref{finiteness} that  $\Ext ^1_R(\ \ , X_2)|_{\pC} \in \modpC$. 
 On the other hand, since  $\Ext ^1_R(\ \ , X_2)$ is a half-exact functor on  $\C$, we see from Lemma \ref{injective} that  $\Ext ^1_R(\ \ , X_2)|_{\pC}$  is an injective object in  $\modpC$. 
\qed
\end{pf}

\begin{itemize}
\item[(Step 3)] For each indecomposable module  $X \in \C$,  the $\pC$-module  $\Ext^1_R ( \ \ , X)|_{\pC}$  is projective in  $\modpC$. 
In particular, $\modpC$  is a Frobenius category. 
\end{itemize}

\begin{pf}
For the proof, we denote the finite set $\ind (\pC)$  by $\{ \underline{W_1}, \ldots , \underline{W_m} \}$ where  $m = |\ind (\pC) |$. 
Setting  $E : = \Ext ^1_R(\ \ , W_i)$  for any one of $i$ ($1 \leq i \leq m$), we want to prove that  $E$  is projective in $\mod\C$. 

Firstly, we show that  $E$  is of finite length as an object in the abelian category  $\modpC$, that is, there is no infinite sequence of strict submodules 
$$
E = E_0 \supset E_1 \supset E_2 \supset \cdots \supset E_n \supset \cdots. 
$$
To show this, set  $W = \oplus _{i=1}^m W_i$  and consider the sequence of $R$-submodules 
$$
E (W) \supset E_1(W) \supset E_2(W) \supset \cdots \supset E_n(W) \supset \cdots$$
 Since we have shown in Lemma \ref{finite length} that 
$E(W) = \Ext ^1_R (W, X)$  is an $R$-module of finite length, this sequence will terminate, i.e. there is an integer  $n$ such that 
$E_n(W) = E_{n+1} (W) = E_{n+2}(W) = \cdots$.
 Since $W$  contains every indecomposable module in  $\pC$, this implies that 
$E_n = E_{n+1} = E_{n+2} = \cdots$  as functors on  $\pC$. 
 Therefore  $E$  is of finite length. 

In particular,  $E$  contains a simple module in $\modpC$  as a submodule. 

Now note that there are only $m$ nonisomorphic indecomposable projective modules in  $\modpC$, in fact they are  
$\pHom _R (\ \ , W_i)|_{\pC} \ (i =1,2,\ldots ,m)$. 
 Corresponding to indecomposable projectives, there are only $m$ nonisomorphic simple modules in $\modpC$ which are 
$$
S_i = \pHom _R (\ \ , W_i)|_{\pC} \ /\ \rad \ \pHom _R (\ \ , W_i)|_{\pC} 
\quad  (i =1,2,\ldots ,m).
$$
 Since we have shown in the steps  1 and 2 that  $\modpC$  is an abelian category with enough projectives and with enough injectives, 
 each simple module  $S_i$  has the injective hull  $I(S_i)$  for  $i=1,2,\ldots ,m$.

Since we have proved that $E$  is an injective module of finite length, 
we see that  $E$  is a finite direct sum of  $I(S_i) \ (i=1,2,\ldots,m)$. 
Since any module in $\modpC$  can be embedded into a direct sum of injective modules of the form  $E = \Ext^1_R(\ \ ,W_i)$, 
we conclude that all nonismorphic indecomposable injective modules in $\modpC$  are  $I(S_i) \ (i =1,2,\ldots ,m)$. 
Note that these are exactly $m$ in number. 

Since  $\modpC$  is a quasi-Frobenius category, any of indecomposable projective modules in $\modpC$  are indecomposable injective. 
 Hence the following two sets coincide: 
$$
\{ \pHom _R (\ \ , W_i)|_{\pC} \ |\ i =1,2,\ldots ,m \} = \{ I(S_i) \ |\ i =1,2,\ldots ,m\}.
$$ 
As a result, every injective module is projective. 
And we have shown that  $\modpC$  is a Frobenius category. 
\qed
\end{pf}

\begin{rem}
We should remark that the proof of the step 3 is the same as the proof of Nakayama's theorem that states the following :

\begin{itemize}
\item[]  Let  $A$  be a finite dimensional algebra over a field. 
 Then $A$  is left selfinjective if and only  $A$  is right selfinjective. 
In particular,  $\mathrm{mod} A$  is a quasi-Frobenius category if and only if  so is $\mathrm{mod} A^{op}$. 
And in this case  $\mathrm{mod} A$  is a Frobenius category. 
\end{itemize}
See \cite{Ya} for example. 
\end{rem}

Now we proceed to the final step of the proof. 
If we prove the following, then the category  $\C$  satisfies all the assumptions in Theorem  \ref{moduloAR}, hence  we obtain  $\C \subseteq \G$. 
And this will complete the proof.

\begin{itemize}
\item[(Step 4)] If  $X \in \C$  such that  $\underline{X} \not\cong \underline{0}$  in  $\pC$, then  we have  $\Ext^1_R(\ \ , X)|_{\pC} \not= 0$. 
\end{itemize}

\begin{pf}
As in the proof of the step 3 we set  
$\ind (\pC) = \{ \underline{W_i} \ | \ i=1,2,\ldots,m\}$. 
It is enough to show that  $\Ext^1_R (\ \ , W_i)|_{\pC} \not= 0$ for each $i$. 
Now assume that  $\Ext^1_R (\ \ , W_1)|_{\pC} =0$  and we shall show a contradiction. 
 In this case, it follows from the step 2 that any module in $\modpC$  can be embedded into a direct sum of copies of  $(m-1)$ modules   $\Ext^1_R (\ \ , W_2)|_{\pC}, \ldots , \Ext^1_R (\ \ , W_m)|_{\pC}$. 
In particular, any indecomposable injective modules appear in these $(m-1)$ modules as direct summands. 
 But we have shown in the proof of the step 3 that there are $m$ indecomposable injective modules  $I(S_i)\ (i=1,2,\ldots,m)$. 
 Hence at least one of  $\Ext^1_R (\ \ , W_2)|_{\pC}, \ldots , \Ext^1_R (\ \ , W_m)|_{\pC}$  contains two different indecomposable injective modules as direct summands.  
Since  $\modpC$  is a Frobenius category,  we see in particular that it is decomposed nontrivially into a direct sum of projective modules in $\modpC$. 
We may assume that   $\Ext^1_R (\ \ , W_2)|_{\pC}$  is decomposed as 
$$
\Ext^1_R (\ \ , W_2)|_{\pC} \cong 
 \pHom _R ( \ \ , Z_1)|_{\pC} \oplus \pHom _R ( \ \ , Z_2)|_{\pC}, 
$$
where  $Z_1 , Z_2 (\not\cong 0) \in \pC$. 
Then it follows from Lemma \ref{main lemma} that  $W_2$  is isomorphic to a direct summand of  $\Omega Z_1\oplus \Omega Z_2$.  
 But since $W_2$  is indecomposable, we may assume that  $W_2$  is isomorphic to a direct summand of  $\Omega Z_1$. 
 Then  $\Ext ^1_R ( \ \ , W_2)|_{\pC}$  is a direct summand of  $\Ext ^1_R ( \ \ , \Omega Z_1)|_{\pC} \cong \pHom _R ( \ \ , Z_1)|_{\pC}$. 
This is a contradiction, because  $\modpC$  is a Krull-Schmidt category by Lemma  \ref{KS thm}. 
\qed
\end{pf}



\begin{thebibliography}{}

\bibitem{A} M.Auslander,
{\it Representation of algebras}{\it I}, Comm. Algebra {\bf 1} (1974), 177 -- 268: {\it II}, Comm. Algebra {\bf 2}(1974), 269 -- 310. 

\bibitem{ABr} M.Auslander and Bridger,
{\it Stable module theory}, 
Mem. Amer. Math. Soc. {\bf 94}, 1969. 

\bibitem{ABu} M.Auslander and R.-O.Buchweitz,
{\it The homological theory of maximal Cohen-Macaulay approximations\/},
Soci\'et\'e Math\'ematique de France, M\'emoire No. 38 (1989), 5--37.

\bibitem{AR} M.Auslander and I.Reiten,
{\it Representation of algebras}{\it III}, Comm. Algebra {\bf 3} (1975), 239 -- 294: {\it IV}, Comm. Algebra {\bf 5}(1977), 443 -- 518. 


\bibitem{Ya} K.Yamagata, 
{\it Frobenius algebras}. 
Handbook of algebra, edited by M.Hazewinkel, 
North-Holland, Amsterdam,  
Vol. 1 (1996), 841--887. 

\bibitem{Y} Y.Yoshino,
{\it Cohen-Macaulay modules over Cohen-Macaulay rings\/},
 London Math. Soc., Lecture Note Series {\bf 146}, Cambridge U.P., 1990.

\end{thebibliography}
\end{document}